\documentclass{article}

\usepackage{amsmath}
\usepackage{url}

\usepackage{amsthm}
\usepackage{amssymb}

\theoremstyle{plain}
\newtheorem{thm}{Theorem}
\newtheorem{lem}[thm]{Lemma}

\theoremstyle{remark}
\newtheorem{rem}{Remark}

\theoremstyle{definition}

\newcommand{\N}{\mathbb N}     

\renewcommand{\l}{\lambda}
\renewcommand{\th}{\theta}

\begin{document}

\title{A probabilistic proof of a binomial identity}
\author{Jonathon Peterson}
\date{}
\maketitle

\begin{abstract}
We give an elementary probabilistic proof of a binomial identity. The proof is obtained by computing the probability of a certain event in two different ways, yielding two different expressions for the same quantity.  
\end{abstract}

The goal of this note is to give a simple (and interesting) probabilistic proof of the binomial identity
\begin{equation}\label{bid}
 \sum_{k=0}^n \binom{n}{k}(-1)^k \frac{\theta}{\theta+k} = \prod_{k=1}^n \frac{k}{\theta+k}, 
\quad \text{for all } \theta > 0 \text{ and all } n \in \N. 
\end{equation}
If one is only concerned with giving a proof of this equality, other proofs than the probabilistic one given below may be more natural. For instance, a proof may be given by identifying the left side of \eqref{bid} as the evaluation of a hypergeometric function $_2F_1\left(-n,\theta;\theta+1|1\right)$
and then applying the Chu-Vandermonde formula \cite[equation (1.2.9)]{grBHS} to obtain the right side of \eqref{bid}.
Another approach would be to use the Rice integral formulas \cite{fsMTFDRI,kNOAS} to equate the left side of \eqref{bid} with a complex contour integral that can be seen to equal the right side of \eqref{bid}. 
These approaches give short proofs of \eqref{bid}, but they both use a good deal of advanced mathematics.
With a bit of work, one can also obtain an elementary proof of \eqref{bid} using only basic properties of the binomial coefficients and mathematical induction. 

The proof of \eqref{bid} given below arose not in a search for a new proof of this identity, but as a result of some independent probability research \ldots\ and a cluttered desk. 
Being unable to find a probability calculation I had done the day before, I sought to repeat the calculation but obtained a different expression for the same quantity. After some initial confusion, I realized that my computations gave a simple proof of the identity \eqref{bid}. 

\section{Probability theory background.}
 
Before giving the probabilistic proof of \eqref{bid}, I will recall some basic facts from probability theory. 
All of the probability needed for this paper can be found in a basic undergraduate probability book such as \cite{rAFCIP}. 
Recall that a random variable $Y$ has an exponential distribution with parameter $\l$ if 
\[
 P(Y\leq y) = \begin{cases} 1-e^{-\l y} & y\geq 0 \\ 0 & y < 0. \end{cases}
\]
The notation $Y\sim \text{Exp}(\l)$ will be used to denote that $Y$ has an exponential distribution with parameter $\l$. 
The following elementary fact about exponential random variables will be used multiple times in the proof of \eqref{bid}. 
\begin{equation}\label{ExpLaplace}
 E[ e^{-\theta Y} ] = \frac{\l}{\l + \theta} \quad \text{for all } \theta > -\l \text{ if } Y\sim \text{Exp}(\l). 
\end{equation}

Another basic tool from probability theory that will be needed is the method of computing probabilities by conditioning. Suppose that $A$ is an event that depends on some random variable $Z$ and some additional randomness. It is sometimes easier to compute the probability of the event $A$ if $Z$ is fixed (that is, by conditioning on $Z$). Then, the probability of the event $A$ is obtained by averaging the conditional probabilities over all values of $Z$: $P(A) = E[ P(A | Z) ]$. 
As an example of this, suppose $Y \sim \text{Exp}(\l)$ and $Z \sim \text{Exp}(\mu)$ are independent. Then, we can compute $P(Y<Z)$ by conditioning on $Y$. Since $P(Y<Z \, | \, Y=y) = e^{-\mu y}$ we see that 
\[
 P( Y < Z ) = E[ P(Y<Z\,|\,Y) ] = E[ e^{-\mu Y}] = \frac{\l}{\l + \mu},
\]
where the last equality follows from \eqref{ExpLaplace}. 

\section{Probabilistic proof.}

Suppose that $X_1,X_2,\ldots,X_n$ are independent Exp(1) random variables, and let $X = \max_{i\leq n} X_i$. Also, let $T \sim \text{Exp}(\theta)$ be independent of the $X_i$ (and thus also independent of $X$). 
The proof of \eqref{bid} given below shows that both sides of \eqref{bid} are equal to the probability $P(X<T)$. 
The two different representations of this probability can both be obtained by conditioning -- the first by conditioning on $X$ and the second by conditioning on $T$.

\subsection{Conditioning on $X$.}
Since $T\sim \text{Exp}(\theta)$, the conditional probability with respect to $X$ is $P(X<T|X=x) = e^{-\theta x}$. 
Therefore,
\begin{equation}\label{Xcond}
 P(X<T) = E[ P(X<T|X)] = E[ e^{-\theta X} ]. 
\end{equation}
The above exponential moment of $X$ can be computed using the following alternate representation of $X$. 
\begin{lem}\label{maxexp}
Suppose that $X_1,X_2,\ldots X_n$ are independent Exp(1) random variables, and that $X = \max_{i\leq n} X_i$. Then, $X$ has the same distribution as $\sum_{k=1}^n Y_k$, where the $Y_k$ are independent random variables with with $Y_k\sim \text{Exp}(k)$.
\end{lem}
\begin{rem}
 The above representation of $X$ as the sum of independent exponential random variables is not standard material in an undergraduate probability course. However, the proof of Lemma \ref{maxexp} below only uses two basic facts about exponential random variables that typically are part of an undergraduate probability course.
\end{rem}
\begin{proof}
 The proof of this lemma relys on two basic facts about exponential distributions. The first is that exponential distributions are memoryless: if $Z\sim \text{Exp}(\l)$ then $P(Z>t+s \, | \, Z > s) = P(Z>t)$. That is, thinking of an exponential random variable as the random amount of time before an event happens, if the event has not occured by time $t$, then the remaining amount of time before the event occurs is still an exponential distribution with the same parameter. 
The second basic fact needed for the proof of Lemma \ref{maxexp} is that if $Z_1, Z_2,\ldots, Z_k$ are independent Exp(1) random variables, then $\min_{i\leq k} Z_i \sim \text{Exp}(k)$. To see this, note that independence implies
\[
 P\left( \min_{i\leq k} Z_i > t \right) = P\left( \bigcap_{i=1}^k \{Z_i > t\} \right)
= \prod_{i=1}^k P(Z_i > t) = e^{-k t}.
\]

Using these two facts, the proof of Lemma \ref{maxexp} is most easily explained in the following way. 
Suppose that there are $n$ lightbulbs in a room that are all turned on at the same time and that $X_i$ is the amount of time until the $i$-th lightbulb fails. 
Then $X = \max_{i\leq n} X_i$ is the amount of time until all of the lightbulbs have failed. 
By the second fact above, the amount of time until one of the lightbulbs burns out is an exponential random variable with parameter $n$. At this time there are still $n-1$ lightbulbs working, and by the first fact above the remaining lifetime of each of these lightbulbs is an Exp(1) random variable.
Thus $X$ has the same distribution as the sum of an Exp($n$) random variable and an independent random variable $X'$ that is the maximum of $n-1$ independent Exp(1) random variables. The conclusion of the lemma then follows by induction on $n$. 
\end{proof}

Applying Lemma \ref{maxexp} to \eqref{Xcond} implies that
\begin{equation}\label{productformula}
 P(X<T) = E\left[ e^{-\theta \sum_{k=1}^n Y_k } \right] = \prod_{k=1}^n E\left[ e^{-\theta Y_k} \right] = \prod_{k=1}^n \frac{k}{k+\theta}, 
\end{equation}
 where the second equality follows from the independence of the $Y_k$ and the last equality follows from \eqref{ExpLaplace}.

\subsection{Conditioning on $T$.}
A different expression for $P(X<T)$ can be obtained by conditioning on $T$ instead. 
Since $\{ X < t \} = \bigcap_{i=1}^n \{ X_i < t \}$ and the $X_i$ are independent, it follows that
\[
 P(X<T \, | \, T=t) = \prod_{i=1}^n P(X_i < t) = \left( 1- e^{-t} \right)^n = \sum_{k=1}^n \binom{n}{k} (-1)^k e^{-kt}.
\]
Then, taking expectations with respect to $T$ gives
\begin{align}
 P(X<T) 
= E\left[\sum_{k=1}^n \binom{n}{k} (-1)^k e^{-kT}\right] 
&= \sum_{k=1}^n \binom{n}{k} (-1)^k E[e^{-kT}] \nonumber \\
&= \sum_{k=1}^n \binom{n}{k} (-1)^k \frac{\theta}{\theta+k}, \label{sumformula}
\end{align}
where the second equality is from the linearity of expected values and the last equality is from \eqref{ExpLaplace}. 
The proof of the binomial identity \eqref{bid} is then completed by combining \eqref{productformula} and \eqref{sumformula}.

\section{Generalizations.}\label{Generalization}

Since this probabilistic proof of \eqref{bid} was constructed quite by accident, it is difficult to use this method to prove a given binomial identity. 
However, the above method can be used to \emph{discover} other interesting binomial identities by making changes to the original probability being computed. 
For instance, let $X$ and $T$ be as above and let $T_2 = T + T'$ where $T'$ is an independent Exp($\th$) random variable. Then, computing $P(X<T_2)$ two different ways will give the following identity
\[
 \sum_{k=0}^n \binom{n}{k} (-1)^k \left( \frac{\th}{\th+k} \right)^2 = \left( \prod_{k=1}^n \frac{k}{\th+k} \right) \left( 1 + \sum_{k=1}^n \frac{\th}{\th+k} \right) .
\]
The left side of the above identity is obtained by first conditioning on $T_2$, while the right side is obtained by conditioning on $X$. 

Even more generally, for any nonnegative integer $m$ one can obtain the identity 
\begin{align*}
& \sum_{k=0}^n \binom{n}{k} (-1)^k \left( \frac{\theta}{\theta+k} \right)^{m} \nonumber \\
& = \left( \prod_{k=1}^n \frac{k}{\th+k} \right) \left( 1 + \sum_{j=1}^{m-1} \, \,  \sum_{1\leq k_1 \leq \ldots \leq k_{j} \leq n}  \frac{\th^{j}}{ (\th+k_1)(\th+k_2)\cdots(\th +k_{j}) } \right)
\end{align*}
by computing $P(X < T_m)$ when $T_m$ is the sum of $m$ independent Exp($\th$) random variables (note that $T_m$ is a Gamma($m$,\ $\th$) random variable). 
Again, these identities could be proved using the Rice integral formulas, but the interested reader may enjoy using the probabilistic method above to prove these identities instead.

\paragraph{Acknowledgments.} 
The probability calculation that led to the above proof was done while doing research supported by National Science Foundation grant DMS-0802942.
I am very grateful to the NSF for this support.


\bigskip

\noindent\textit{Department of Mathematics,
Purdue University, 150 N.\ University Street, West Lafayette, IN 47907\\
peterson@math.purdue.edu}

\end{document}